 \documentclass[preprint,12pt]{elsarticle}

\usepackage{amssymb, amsfonts, amsmath, graphicx, algorithm, algorithmic}

\newtheorem{thm}{Theorem}
\newtheorem{cor}[thm]{Corollary}
\newtheorem{lem}[thm]{Lemma}

\newtheorem{prop}[thm]{Proposition}

\newdefinition{defn}{Definition}

\newdefinition{rem}{Remark}

\newproof{pf}{Proof}

\begin{document}

\begin{frontmatter}



\title{Weighted matrix means and symmetrization procedures\tnoteref{label1}}

\tnotetext[label1]{I would like to thank prof. Yongdo Lim for inviting me as a visitor to Kyungpook National University, where some parts of this research has been carried out.}

\author[label2]{Mikl\'os P\'alfia}
\ead{palfia.miklos@aut.bme.hu}

\date{\today}

\address[label2]{Department of Automation and Applied Informatics, Budapest University of Technology and Economics, H-1521 Budapest, Hungary}

\begin{abstract}
Here we prove the convergence of the Ando-Li-Mathias and Bini-Meini-Poloni procedures for matrix means. Actually it is proved here that for a two-variable function which maps pairs of positive definite matrices to a positive definite matrix and is not greater than the square mean of two positive definite matrices, the Ando-Li-Mathias and Bini-Meini-Poloni procedure converges. In order to be able to set up the Bini-Meini-Poloni procedure, a weighted two-variable matrix mean is also needed. Therefore a definition of a two-variable weighted matrix mean corresponding to every symmetric matrix mean is also given. It is also shown here that most of the properties considered by Ando, Li and Mathias for the $n$-variable geometric mean hold for all of these $n$-variable maps that we obtain by this two limiting process for all two-variable matrix means. As a consequence it also follows that the Bini-Meini-Poloni procedure converges cubically for every matrix mean.
\end{abstract}

\begin{keyword}
means \sep operator means \sep geometric mean
\PACS 15A24 \sep 15A45 \sep 26E60 \sep 47A64
\end{keyword}
\end{frontmatter}

\section{Introduction}
Recently several researchers were considering the extension of means to multiple variables. Most of the contributions aimed at the geometric mean of positive definite matrices. Every extension is based on the two-variable form of the geometric mean defined by the theory of Kubo and Ando which characterizes a wide family of two-variable matrix means. These two-variable means are isomorphic to normalized operator monotone functions, therefore they are well characterized by the theory.

The geometric mean is of great interest, since it is less trivial than the arithmetic or harmonic mean of positive matrices and at the same time it also has a corresponding nonpositively curved Riemannian metric for which it is the midpoint operation. This metric provides a great tool to extend the two-variable geometric mean of matrices to several variables. For instance with the aid of this corresponding metric structure Ando-Li-Mathias \cite{ando}, Bini-Meini-Poloni \cite{bini}, Jung-Lee-Yamazaki \cite{jung} and Moakher \cite{moakher} were able to define the geometric mean of several positive operators. The first three constructions are based on different iterative procedures which define the mean as a limit point for $n\geq 3$ variables. The geometric $n$-mean of Moakher or the so called Riemannian mean is defined as the center of mass of the $n$ matrices.

All of the procedures extend to metric spaces where the two-mean fulfills some geometrical properties, or the space itself has some additional geometrical properties. For example Lawson and Lim in \cite{lawsonlim} were able to extend the Ando-Li-Mathias procedure to metric spaces where the two-mean fulfills some geometric property expressed in an inequality closely related to Busemann nonpositive curvature. Actually the Ando-Li-Mathias procedure was considered by Petz and Temesi \cite{petz} as extensions of matrix means to several variables, however the convergence of the procedure was only proved there for orderable tuples.

Recently the author in \cite{palfia} defined means of several matrices based on an iterative procedure relying on a properly chosen infinite sequence of graphs $G$. This procedure extends every two-variable matrix mean given by the Kubo-Ando theory to several variables. The applicability of the extension has been proved in full generality. This procedure was later extended to geodesic metric spaces or more specifically $C_k$-domains where the two-mean is considered as the midpoint map of the space \cite{palfia2}. In this article we prove that the Ando-Li-Mathias (ALM) and Bini-Meini-Poloni (BMP) procedure converges for any matrix mean. Since the BMP-procedure requires the existence of a weighted matrix mean, we therefore give a sufficient definition of a two-variable weighted matrix mean which agrees with the theory of Kubo-Ando.

\section{Matrix means in two-variables}
In this section we recall the Kubo-Ando theory of matrix means. We investigate some properties of these means that will be useful in the following sections. By the foundational article \cite{kubo} we consider two-variable functions on the open convex cone of hermitian positive-definite $r\times r$ matrices $P(r)$.
\begin{defn}[Matrix mean]\label{mean}
Let $M:P(r)\times P(r)\mapsto P(r)$. Then $M(A,B)$ is called a matrix mean if the following conditions hold
\begin{enumerate}
\renewcommand{\labelenumi}{(\roman{enumi})}
\item $M(I,I)=I$ where $I$ denotes the identity,
\item if $A\leq A'$ and $B\leq B'$, then $M(A,B)\leq M(A',B')$,
\item $CM(A,B)C\leq M(CAC,CBC)$ for all hermitian $C$,
\item if $A_n\downarrow A$ and $B_n\downarrow B$ then $M(A_n,B_n)\downarrow M(A,B)$.
\newcounter{enumi_saved}
\setcounter{enumi_saved}{\value{enumi}}
\end{enumerate}
\end{defn}
As consequences of the above properties we can derive further properties of matrix means:
\begin{enumerate}
\renewcommand{\labelenumi}{(\roman{enumi})}
\setcounter{enumi}{\value{enumi_saved}}
\item if $A\leq B$ then $A\leq M(A,B)\leq B$,
\item $SM(A,B)S^{*}\leq M(SAS^{*},SBS^{*})$ for $S$ not necessarily hermitian,
\item $M(CAC^{*},CBC^{*})=CM(A,B)C^{*}$ if $C$ is invertible,
\item $M(A,B)+M(C,D)\leq M(A+C,B+D)$,
\setcounter{enumi_saved}{\value{enumi}}
\end{enumerate}
and the most important consequence of the Kubo-Ando theory of means is property
\begin{enumerate}
\renewcommand{\labelenumi}{(\roman{enumi})}
\setcounter{enumi}{\value{enumi_saved}}
\item $M(A,B)=A^{1/2}f\left(A^{-1/2}BA^{-1/2}\right)A^{1/2}$,
\end{enumerate}
where $f$ is a non-negative operator monotone function which is normalized so $f(I)=I$. The last property provides an isomorphism between means of two hermitian positive definite matrices and normalized operator monotone functions. As an immediate consequence we have that every matrix mean $M(A,B)$ is continuous. The proof of properties (vi)-(ix) can be found in \cite{kubo}. Property (v) is a direct consequence of property (i),(ii) and (ix), since if $A\leq B$ then
\begin{equation}\label{inbetween}
A=M(A,A)\leq M(A,B)\leq M(B,B)=B\text{.}
\end{equation}

Kubo-Ando theory considers not necessarily symmetric matix means, so generally $M(A,B)\neq M(B,A)$. Although it turns out, according to \cite{kubo}, that if $M(A,B)$ is symmetric then we have
\begin{equation}\label{symmetricbounds}
\left(\frac{A^{-1}+B^{-1}}{2}\right)^{-1}\leq M(A,B)\leq \frac{A+B}{2}\text{.}
\end{equation}
In other words the smallest symmetric mean is the harmonic mean, the largest is the arithmetic mean with respect to the standard positive definite order. This fact is the consequence of the isomorphism between normalized operator monotone functions and matrix means.

\section{Weighted matrix means in two-variables}
Although Kubo-Ando theory characterizes all possible matrix means, it does not clarify how weighted means correspond to symmetric ones. In some special cases we have natural correspondence between a symmetric mean and some other weighted means. In most cases this correspondence is based on some affine geometric structures. For instance the weighted arithmetic mean
\begin{equation}\label{weightedharmonic}
A_t(A,B)=(1-t)A+tB
\end{equation}
is the geodesic line connecting $A$ and $B$ with respect to the standard Euclidean metric corresponding to the space of complex squared matrices. This metric has the form $\left\langle X,Y\right\rangle_p=Tr\left\{XY^{*}\right\}$ on each tangent space. Another related example is the weighted harmonic mean
\begin{equation}\label{weightedarithmetic}
H_t(A,B)=\left((1-t)A^{-1}+tB^{-1}\right)^{-1}
\end{equation}
and it is also a geodesic with respect to the metric $\left\langle X,Y\right\rangle_p=Tr\left\{p^{-2}Xp^{-2}Y\right\}$, where $p\in P(r)$ and $X,Y$ are hermitian.

There is another important example which is the weighted geometric mean
\begin{equation}
G_t(A,B)=A^{1/2}\left(A^{-1/2}BA^{-1/2}\right)^tA^{1/2}\text{.}
\end{equation}
The corresponding Riemannian metric is $\left\langle X,Y\right\rangle_p=Tr\left\{p^{-1}Xp^{-1}Y\right\}$, where again $p\in P(r)$ and $X,Y$ are hermitian. This manifold is nonpositively curved while the other two is Euclidean, therefore they have zero curvature.

In the above three cases the structure of the underlying affinely connected manifold gives us the natural correspondence between the symmetric and weighted counterparts of the matrix mean. In these cases the geodesic lines $\gamma:[0,1]\mapsto P(r)$ provides us the weighted matrix means. We are going to use the following definition to identify these special means.
\begin{defn}[Affine mean]\label{affinemean}
An affine mean $M:W^2\mapsto W$ is a geodesic midpoint operation $M(A,B)=\exp_A(1/2\log_A(B))$ on a smooth manifold $W$ equipped with an affine connection, where $B$ is assumed to be in the injectivity radius of the exponential map $\exp_A(x)$ of the connection given at the point $A$. The mapping $\log_A(x)$ is just the inverse of the exponential map at the point $A\in W$.
\end{defn}
In full generality it is not known whether such affine manifolds exists for every symmetric matrix mean where the mean is the midpoint operation on the manifold. This appears to be a delicate problem however we can borrow some ideas based on the above three examples to define weighted mean counterparts for every symmetric matrix mean.

First of all we define a procedure for every symmetric matrix mean $M(A,B)$ and for all $t\in[0,1]$. Our procedure will be based on the fact that every $t\in[0,1]$ can be approximated by dyadic rationals $\frac{m}{2^n}$ since dyadic rationals are dense in $[0,1]$.
\begin{defn}[Weighted mean process]\label{weightedmeanprocess}
Let $M(\cdot,\cdot)$ be a symmetric matrix mean, $A,B\in P(r)$ and $t\in[0,1]$. Let $a_0=0$ and $b_0=1$, $A_0=A$ and $B_0=B$. Define $a_{n}, b_{n}$ and $A_{n}, B_{n}$ recursively by the following procedure for all $n=0,1,2,\dots$:
\begin{algorithmic}
\IF{$a_n=t$}
	\STATE $a_{n+1}=a_n$ and $b_{n+1}=a_n$, $A_{n+1}=A_n$ and $B_{n+1}=A_n$
\ELSIF{$b_n=t$}
	\STATE $a_{n+1}=b_n$ and $b_{n+1}=b_n$, $A_{n+1}=B_n$ and $B_{n+1}=B_n$
\ELSIF{$\frac{a_n+b_n}{2}\leq t$}
	\STATE $a_{n+1}=\frac{a_n+b_n}{2}$ and $b_{n+1}=b_n$, $A_{n+1}=M(A_n,B_n)$ and $B_{n+1}=B_n$
\ELSE
	\STATE $b_{n+1}=\frac{a_n+b_n}{2}$ and $a_{n+1}=a_n$, $B_{n+1}=M(A_n,B_n)$ and $A_{n+1}=A_n$
\ENDIF
\end{algorithmic}
According to the above $a_{n+1}, b_{n+1}$ and $A_{n+1}, B_{n+1}$ are clearly defined with respect to $a_{n}, b_{n}$ and $A_{n}, B_{n}$ recursively.
\end{defn}
This algorithm may also be regarded as a kind of binary search with recurrence relation:
\begin{algorithmic}
\IF{$t=\frac{t_1+t_2}{2}$}
	\STATE $M_t(A,B)=M\left(M_{t_1}(A,B),M_{t_2}(A,B)\right)$
\ENDIF
\end{algorithmic}
\begin{thm}\label{weightedmeanprocessconv}
The sequences $A_n$ and $B_n$ given in Definition \ref{weightedmeanprocess} are convergent and have the same limit point.
\end{thm}
\begin{pf}
In the case if $t=m2^{-k}$ for some integer $m$ and $k$ then there is nothing to prove, the procedure converges after finite steps. So suppose that $t$ is not a dyadic rational. We will make use of the following multiplicative metric on $P(r)$ \cite{ando}
\begin{equation}\label{defrab}
R(A,B)=\max\left\{\rho(A^{-1}B),\rho(B^{-1}A)\right\}
\end{equation}
for all $A,B\in P(r)$ and $\rho(A)$ denotes the spectral radius of $A$. The above metric has the following properties \cite{bini}
\begin{enumerate}
\item $R(A,B)\geq 1$,
\item $R(A,B)=1$ iff $A=B$,
\item $R(A,C)\leq R(A,B)R(B,C)$,
\item $R(A,B)^{-1}A\leq B\leq R(A,B)A$,
\item $\left\|A-B\right\|\leq (R(A,B)-1)\left\|A\right\|$.
\end{enumerate}
Since $R(A,B)=R(I,A^{-1/2}BA^{-1/2})$ we have
\begin{equation}
R\left(A,M(A,B)\right)=R\left(I,f\left(A^{-1/2}BA^{-1/2}\right)\right)\text{,}
\end{equation}
where $f(t)$ is the corresponding normalized operator monotone function. Now since $M(A,B)$ is symmetric \eqref{symmetricbounds} holds. From this for the corresponding normalized operator monotone function $f(t)$ we have
\begin{equation}
2\left(I+X^{-1}\right)^{-1}\leq f(X)\leq \frac{I+X}{2}\text{.}
\end{equation}
This also yields that $R(I,f(X))\leq \max\left\{\rho(\frac{I+X}{2}),\rho(\frac{I+X^{-1}}{2})\right\}=\frac{1+R(I,X)}{2}$ for every $X\in P(r)$, so
\begin{equation}
R\left(A,M(A,B)\right)\leq \frac{1+R(A,B)}{2}\text{.}
\end{equation}
By the above inequality we can easily conclude the following for the sequences $A_n, B_n$
\begin{equation}
\begin{split}
R(A_{n+1},B_{n+1})&\leq \frac{1+R(A_n,B_n)}{2}=1+\frac{1}{2}\left[R(A_n,B_n)-1\right]\\
R(A_n,B_n)&\leq 1+\frac{1}{2^n}\left[R(A_0,B_0)-1\right]\\
R(A_n,A_{n+1})&\leq 1+\frac{1}{2}\left[R(A_n,B_n)-1\right]\\
R(A_n,A_{n+1})&\leq 1+\frac{1}{2^n}\left[R(A_0,B_0)-1\right]\text{.}
\end{split}
\end{equation}
There exists $K\in P(r)$ such that $A\leq K, B\leq K$ and by property (ii) of matrix means $A_n\leq K, B_n\leq K$ so by property 5. of $R(\cdot,\cdot)$
\begin{equation}
\begin{split}
\left\|A_{n+1}-A_{n}\right\|&\leq (R(A_{n+1},A_{n})-1)\left\|K\right\|\\
\left\|A_{n+1}-A_{n}\right\|&\leq \frac{1}{2^n}\left[R(A_0,B_0)-1\right]\left\|K\right\|\\
\sum_{n=0}^{\infty}\left\|A_{n+1}-A_{n}\right\|&\leq \sum_{n=0}^{\infty}\frac{1}{2^n}\left[R(A_0,B_0)-1\right]\left\|K\right\|=\\
&=2\left[R(A_0,B_0)-1\right]\left\|K\right\|\text{.}
\end{split}
\end{equation}
This means that $A_n$ is a Cauchy sequence therefore convergent and by the above we also have that $\left\|A_{n}-B_{n}\right\|\to 0$ so both $A_n$ and $B_n$ converge to the same limit point.$\square$
\end{pf}
We will base our weighted mean on the above theorem.
\begin{defn}[Weighted mean]
The common limit point of $A_n,B_n$ in Theorem \ref{weightedmeanprocessconv} will be denoted by $M_t(A,B)$ and from now on in the article is considered as the corresponding weighted mean to a symmetric matrix mean $M(\cdot,\cdot)$.
\end{defn}
What are the properties of this weighted mean? First of all it is not hard to prove the following 
\begin{prop}\label{correctcorrespondings}
$M_t(A,B)$ yields the correct corresponding weighted means in the case of the arithmetic, geometric, harmonic means.
\end{prop}
The above is a consequence of the affine geodesy of the corresponding manifolds mentioned above. There are further important properties which are fulfilled by $M_t(A,B)$:
\begin{prop}\label{weightedmeanproperties}
$M_t(A,B)$ for $A,B\in P(r)$ and $t\in[0,1]$ fulfills the following properties
\begin{enumerate}
\renewcommand{\labelenumi}{(\roman{enumi}')}
\item $M_t(I,I)=I$,
\item if $A\leq A'$ and $B\leq B'$, then $M_t(A,B)\leq M_t(A',B')$,
\item $CM_t(A,B)C\leq M_t(CAC,CBC)$,
\item if $A_n\downarrow A$ and $B_n\downarrow B$ then $M_t(A_n,B_n)\downarrow M_t(A,B)$,
\item if $N(A,B)\leq M(A,B)$ then $N_t(A,B)\leq M_t(A,B)$,
\item $M_{1/2}(A,B)=M(A,B)$,
\item $M_t(A,B)$ is continuous in $t$,
\end{enumerate}
\end{prop}
\begin{pf}
Property (i') and (ii') are trivial consequences of the similar properties for symmetric matrix means.

For property (iii') consider $A'=CAC$ and $B'=CBC$ and start the procedure in the definition of $M_t(\cdot,\cdot)$ for the pair $A,B$ and $A',B'$. Then we have $CA_{1}C=CM(A_0,B_0)C\leq M(CA_0C,CB_0C)=A'_1$ if $t>1/2$ or $CB_{1}C=CM(A_0,B_0)C\leq M(CA_0C,CB_0C)=B'_1$. Now for every $n$ we use property (ii) for symmetric matrix means so we have $CA_nC \leq A'_n$ and $CB_nC \leq B'_n$ for every $n\geq 1$. Taking the limits we conclude the assertion of property (iii').

What immediately follows from this property is that $M(CXC^{*},CYC^{*})=CM(X,Y)C^{*}$ for all invertible $C$. We will use this to show that $M_t(A,B)$ is continuous in $A,B$ so by property (ii') and the continuity we get property (iv') as a consequence. Actually we prove more, we will show that for a function
\begin{lem}\label{monotonecont}
$F:P(r)^n\mapsto P(r)$ which satisfies properties
\begin{enumerate}
\item if $X_i\leq X'_i$ for all $i$, then $F(X_1,\dots, X_n)\leq F(X'_1,\dots, X'_n)$,
\item $F(cX_1,\dots, cX_n)=cF(X_1,\dots, X_n)$ for real $c>0$,
\end{enumerate}
it follows that $F$ is continuous.
\end{lem}
\begin{pf}
We know that for a function $f:Y_1\to Y_2$ between two metric spaces $(Y_1,d_1)$ and $(Y_2,d_2)$ sequential continuity and the usual topological continuity are equivalent. A proof can be found for example in \cite{maddox}. We will show that sequential continuity holds therefore arriving at the desired result.

We extend the metric $R(\cdot,\cdot)$ to $P(r)^n$ as follows. Let $X=\left(X_1,\dots,X_n\right)\in P(r)^n$ and $Y=\left(Y_1,\dots,Y_n\right)\in P(r)^n$, then we define
\begin{equation}
R_n(X,Y)=\max_{1\leq i\leq n}\left\{R(X_i,Y_i)\right\}\text{.}
\end{equation}
Now suppose we have a convergent sequence of tuples $X^{k}=\left(X_1^{k},\ldots,X_n^{k}\right)\in P(r)^n$ for which $\left(X_1^{k},\ldots,X_n^{k}\right)\to \left(X_1,\ldots,X_n\right)=X\in P(r)^n$. Using property (iv) of $R(\cdot,\cdot)$ we have the following inequalities
\begin{equation}
R_n(X^{k},X)^{-1}X_i^{k}\leq X_i\leq R_n(X^{k},X)X_i^{k}\text{.}
\end{equation}
Now applying the monotonicity property 1. of $F$, we have with the notation $c_k=R_n(X^{k},X)$ the following
\begin{equation}
F\left(c_k^{-1}X_1^{k},\ldots,c_k^{-1}X_n^{k}\right)\leq F\left(X_1,\ldots,X_n\right)\leq F\left(c_kX_1^{k},\ldots,c_kX_n^{k}\right)\text{.}
\end{equation}
Using property 2. we conclude that
\begin{equation}
c_k^{-1}F\left(X_1^{k},\ldots,X_n^{k}\right)\leq F\left(X_1,\ldots,X_n\right)\leq c_kF\left(X_1^{k},\ldots,X_n^{k}\right)\text{.}
\end{equation}
Taking the limit $k\to \infty$ we have $c_k\to 1$. This shows that
\begin{equation}
\lim_{k\to\infty}F\left(X_1^{k},\ldots,X_n^{k}\right)=F\left(\lim_{k\to\infty}X_1^{k},\ldots,\lim_{k\to\infty}X_n^{k}\right)
\end{equation}
which is sequential continuity therefore arriving at the assertion of the lemma.
\end{pf}

Since the properties needed in the above lemma holds for $M_t(A,B)$ we conclude by the above that property (iv') must also hold. At this point we already have by the Kubo-Ando theory of matrix means that $M_t(A,B)$ is a matrix mean as well, so it fulfills the additional properties (v)-(ix). Consequently it has a representation with a normalized operator monotone function.

Property (v') is an easy consequence of repeated usage of property (ii) for matrix means for every $n$. Property (vi') is also trivial.

To prove property (vii') we have to do a bit more work. We have to show that if $\left|t_1-t_2\right|\to 0$ then also $\left\|M_{t_1}(A,B)-M_{t_2}(A,B)\right\|\to 0$. Suppose $t_1<t_2$ and take the smallest $j$ for which we have $t_1\leq m2^{-j}\leq t_2$ for some $m$. Let us set up the iterative procedure in Definition \ref{weightedmeanprocess} on $A,B$ with $t_1$ and $t_2$ respectively. Let us denote the yielded matrix sequences in the case of $t_1$ with $A^{t_1}_i,B^{t_1}_i$ and in the case of $t_2$ with $A^{t_2}_i,B^{t_2}_i$ and similarly for the numbers with $a^{t_1}_i,b^{t_1}_i$ and $a^{t_2}_i,b^{t_2}_i$. Notice that the iterative procedure in the $j$th step for $t_1$ will yield $b^{t_1}_j=m2^{-j}$ and similarly $a^{t_2}_j=m2^{-j}$ in the $j$th step for $t_2$. Suppose $t_1\neq m2^{-j}$. Then there exists $i\geq j$ such that $a^{t_1}_i\leq t_1\leq b^{t_1}_i$ but $(a^{t_1}_i+b^{t_1}_i)/2\geq t_1$, this means that $b^{t_1}_{i+1}\neq b^{t_1}_{i}$. If $t_1=m2^{-j}$ then we have $a^{t_1}_p=b^{t_1}_p=t_1$ for $p>j$ and in this case we define $i:=+\infty$. Similarly either there exists a smallest $l\geq j$ such that $a^{t_2}_{l+1}\neq a^{t_2}_{l}$, or we have $t_2=m2^{-j}$ so $a^{t_2}_p=b^{t_2}_p=t_2$ for $p>j$ and again then we define $l:=+\infty$. Notice that $i$ and $l$ cannot be infinite at the same time, so we define $k:=\min\left\{i,l\right\}$. It is easy to see that as $t_1\to t_2$, $k\to \infty$. We also have that $B^{t_1}_k=A^{t_2}_k$, so we can bound the distance of the limit points $M_{t_1}(A,B)$ and $M_{t_2}(A,B)$ from $B^{t_1}_k=A^{t_2}_k$ as follows:
\begin{equation}
\begin{split}
\left\|B^{t_1}_k-\lim_{j\to\infty} B^{t_1}_j\right\|&\leq \sum_{i=k}^{\infty}\left\|B^{t_1}_{i+1}-B^{t_1}_i\right\|\\
\sum_{i=k}^{\infty}\left\|B^{t_1}_{i+1}-B^{t_1}_i\right\|&\leq \frac{1}{2^k}2\left[R(A_0,B_0)-1\right]\left\|K\right\|=\\
&=\frac{1}{2^{k-1}}\left[R(A_0,B_0)-1\right]\left\|K\right\|\text{.}
\end{split}
\end{equation}
We also have the same bound for $\left\|A^{t_2}_k-\lim_{j\to\infty} A^{t_2}_j\right\|$. Since $B^{t_1}_k=A^{t_2}_k$, we have
\begin{equation}
\begin{split}
\left\|M_{t_1}(A_0,B_0)-M_{t_2}(A_0,B_0)\right\|\leq &\left\|M_{t_1}(A_0,B_0)-B^{t_1}_k\right\|+\\
&+\left\|M_{t_2}(A_0,B_0)-A^{t_2}_k\right\|\leq\\
\leq &\frac{1}{2^{k-2}}\left[R(A_0,B_0)-1\right]\left\|K\right\|\text{.}
\end{split}
\end{equation}
Since $k\to \infty$ as $t_1\to t_2$, by the above $\left\|M_{t_1}(A,B)-M_{t_2}(A,B)\right\|\to 0$. $\square$
\end{pf}

By the above proposition we have that $M_t(A,B)$ is a continuous function in $t$. So $M_t(A,B)$ is a one parameter family of matrix means corresponding to every symmetric matrix mean. Since every matrix mean by virtue of property (ix) is representable by a normalized operator monotone function $f(x)$, we may represent such one parameter family of matrix means by a one parameter family of normalized operator monotone functions $f_t(x), t\in [0,1]$. So in other words we have the following
\begin{cor}
For every symmetric matrix mean $M(A,B)$ there is a corresponding one parameter family of weighted means $M_t(A,B)$ for $t\in [0,1]$. Let $f(x)$ be the normalized operator monotone function corresponding to $M(A,B)$. Then similarly we have a one parameter family of normalized operator monotone functions $f_t(x)$ corresponding to $M_t(A,B)$. The family $f_t(x)$ is continuous in $t$, and $f_0(x)=1$ and $f_1(x)=x$ are the two extremal points, so $f_t(x)$ interpolates between these two points.
\end{cor}


Based on this phenomenon we can conclude the following
\begin{prop}\label{arithmeticmatrixmean}
Let $M(A,B)$ be a symmetric matrix mean. Then
\begin{equation}
\left((1-t)A^{-1}+tB^{-1}\right)^{-1}\leq M_t(A,B)\leq (1-t)A+tB\text{,}
\end{equation}
where $M_t(A,B)$ is the weighted version of $M(A,B)$.
\end{prop}
\begin{pf}
For every symmetric matrix mean $M(A,B)$ we have by \eqref{symmetricbounds} that
\begin{equation}
\left(\frac{A^{-1}+B^{-1}}{2}\right)^{-1}\leq M(A,B)\leq \frac{A+B}{2}\text{.}
\end{equation}
Now what follows from Propostion~\ref{correctcorrespondings} is that the harmonic mean on the left hand side above has the weighted harmonic mean $H_t(A,B)$ defined in \eqref{weightedharmonic} as its weighted counterpart, and similarly we have the weighted arithmetic mean $A_t(A,B)$ defined in \eqref{weightedarithmetic} as the weighted counterpart for the arithmetic mean on the right hand side above. Thus by property (v') in Proposition \ref{weightedmeanproperties} and the above inequality we have
\begin{equation}
H_t(A,B)\leq M_t(A,B)\leq A_t(A,B)\text{.}
\end{equation}
$\square$
\end{pf}

We are going to use the above definition $M_t(A,B)$ of a weighted matrix mean to set up the Bini-Meini-Poloni procedure for every symmetric matrix mean, but before we do that we turn to the Ando-Li-Mathias procedure in the following section and prove its convergence for every symmetric matrix mean.

\section{Ando-Li-Mathias procedure for every matrix mean}
In this section we will prove the convergence of the Ando-Li-Mathias procedure for every possible symmetric matrix mean. In order to do that we will generalize the argument given in \cite{palfia} by applying induction. First of all let us recall the Ando-Li-Mathias procedure \cite{ando}:

\begin{defn}\label{almit}[ALM iteration]
Let $X=(X_1^0,\dots,X_n^0)$ where $X_i^0 \in P(r)$ and define the mapping $M(X_1,\dots,X_{n})$ inductively as follows. If $n=2$ assume that $M(X_1,X_2)$ is already given. For general $n>2$ assume that $M(X_1,\dots,X_{n-1})$ is already defined. Then using $M(X_1,\dots,X_{n-1})$, set up the iteration
\begin{equation}
X_i^{l+1}=M\left(Z_{\neq i}\left(X_1^l,\dots,X_n^l\right)\right)\text{,}
\end{equation}
where $Z_{\neq i}(X_1^l,\dots,X_n^l)=X_1^l,\dots,X_{i-1}^l, X_{i+1}^l,\dots,X_n^l$. If the sequences $X_i^{l}$ converge to a common limit point for every $i$, then define
\begin{equation}
\lim_{l\to\infty}X_i^{l}=M(X_1^0,\dots,X_n^0)\text{.}
\end{equation}
\end{defn}

\begin{thm}\label{almthm}
Let $F:P(r)^2\mapsto P(r)$ and suppose that $F(A,B)$ fulfills one of the inequalities below:
\begin{equation}\label{smaller}
\left(\frac{A^{-1}+B^{-1}}{2}\right)^{-1}\leq F(A,B)\leq \left[\frac{A^2+B^2}{2}-\frac{k}{8}(A-B)^2\right]^{1/2}
\end{equation}
for a $k\in(0,2]$, or
\begin{equation}\label{arithmsmaller}
F(A,B)\leq \frac{A+B}{2}\text{.}
\end{equation}
Then in Definition~\ref{almit} starting with $M(A,B):=F(A,B)$, $M(X_1,\dots, X_{n})$ exists for all $n$, in other words the sequences converge to a common limit point for all $n$.
\end{thm}
Before we prove the above theorem we mention a few remarks and several lemmas which we will make use of later. First of all condition \eqref{smaller} might seem a bit strange at first glance although it immediately becomes straightforward if we consider $k=2$, since in this case the right hand side becomes the arithmetic mean. If $k=0$ then the right hand side is the square mean $\left(\frac{A^2+B^2}{2}\right)^{1/2}$. This literally means that the above theorem automatically covers every symmetric matrix mean due to \eqref{symmetricbounds} as a special case.

Now we have to study some properties of the square mean in order to prepare the necessary steps for the proof of the above theorem. First of all one should notice that the square mean is an affine mean. The underlying manifold is a Riemannian manifold defined as a pullback metric of the Euclidean metric $\left\langle A,B\right\rangle_E=Tr\left\{A^*B\right\}$ over the space of squared complex matrices. This Euclidean metric has corresponding distance function
\begin{equation}\label{frobeniusdistance}
\begin{split}
d_E(A,B)^2=\left\langle A-B,A-B\right\rangle_E=\\
=Tr\left\{\left(A-B\right)^{*}\left(A-B\right)\right\}\text{.}
\end{split}
\end{equation}

The isometry is $f(x)=x^2$ and it embeds $P(r)$ into $P(r)$. The distance function of the pullback metric on $P(r)$ is 
\begin{equation}\label{squreddistancefunc}
\begin{split}
d_{1/2}(A,B)^2=\left\langle f(A)-f(B),f(A)-f(B)\right\rangle_E=\\
=Tr\left\{\left(A^2-B^2\right)^{*}\left(A^2-B^2\right)\right\}\text{.}
\end{split}
\end{equation}
The geodesics of this metric are of the form
\begin{equation}
\gamma_{A,B}(t)=f^{-1}\left[(1-t)f(A)+tf(B)\right]=\left[(1-t)A^2+tB^2\right]^{1/2}\text{.}
\end{equation}
This shows that square mean is an affine mean, so the weighted mean process $M_t(A,B)$ for the square mean gives back the corresponding point on the geodesic above. Furthermore since the above metric is a pullback of a Euclidean metric, it is also Euclidean.

Actually the isometry $f(x)$ can be chosen arbitrarily, particularly any diffeomorphism will suffice. We are going to derive some properties of the ALM-procedure with $F(A,B)=\left(\frac{A^2+B^2}{2}\right)^{1/2}$ and $F(A,B)=\left(\frac{A^{-1}+B^{-1}}{2}\right)^{-1}$ on $P(r)$ endowed with the above corresponding pullback metrics. We are going to denote the general pullback of the distance function $d_E(\cdot,\cdot)$ for an arbitrary $f$ by
\begin{equation}\label{pullbackdistance}
\begin{split}
d_{f}(A,B)^2=\left\langle f(A)-f(B),f(A)-f(B)\right\rangle_E=\\
=Tr\left\{\left[f(A)-f(B)\right]^{*}\left[f(A)-f(B)\right]\right\}\text{.}
\end{split}
\end{equation}

The metric space $P(r)$ with the distance function \eqref{pullbackdistance} is Euclidean, since its metric is a pullback metric of the standard Euclidean metric on the space of complex $r\times r$ matrices. Let $x_i\in P(r)$ for $i\in \left\{1,\dots,n\right\}$ and define $S=\left\{x_1,\dots,x_n\right\}$. Then according to \cite{karcher} the function
\begin{equation}
b(x)=\sum_{i=1}^nd(x,x_i)^2
\end{equation}
has a minimum for $d(\cdot,\cdot)=d_{f}(\cdot,\cdot)$ and this minimal value is attained at a unique point $\hat{x}$ which is called the Riemann centroid of $S$. Moreover the centroid can be explicitly given for these metric spaces on $P(r)$.

\begin{prop}\label{centroidsquare}
The unique minimizer $\hat{x}$ of the function
\begin{equation}
b(x)=\sum_{i=1}^nd(x,x_i)^2
\end{equation}
for the distance function \eqref{pullbackdistance} is given as
\begin{equation}
\hat{x}=f^{-1}\left(\frac{\sum_{i=1}^nf(x_i)}{n}\right)\text{.}
\end{equation}
\end{prop}
\begin{pf}
Since the corresponding metric is a pullback of the Euclidean metric over the space of squared complex matrices it is also Euclidean. Using the isometric embedding $f(x)$, the object function of the minimization problem is of the form
\begin{equation}
\sum_{i=1}^nd_{f}(x,x_i)^2=\sum_{i=1}^nd_E\left(f(x),f(x_i)\right)^2\text{.}
\end{equation}
But since the Riemann centroid of the set $S=\left\{f(x_1),\dots,f(x_n)\right\}$ in the Euclidean space of squared complex matrices is the arithmetic mean of the points $\left\{f(x_1),\dots,f(x_n)\right\}$, therefore
\begin{equation}
a=\frac{\sum_{i=1}^nf(x_i)}{n}
\end{equation}
minimizes the functional $\sum_{i=1}^nd_E\left(x,f(x_i)\right)^2$, so $\hat{x}=f^{-1}(a)$ minimizes $\\ \sum_{i=1}^nd_E\left(f(x),f(x_i)\right)^2$.
$\square$
\end{pf}

If we perform one ALM-iteration step on $n$ points in the space $P(r)$ with this centroid map then the iteration leaves the Riemann centroid of the points invariant.

\begin{prop}\label{invariantsquaremean}
Let $X_i^0\in P(r)$ for $i=1,\ldots,n$. Then the ALM-procedure (Definition~\ref{almit}) set up on the matrices $X_1^0,\dots,X_n^0$ with the $n-1$ variable function $M(x_1,\dots,x_{n-1})=f^{-1}\left(\frac{\sum_{i=1}^{n-1}f(x_i)}{n-1}\right)$ leaves the Riemann centroid of the points $X_1^0,\dots,X_n^0$ invariant with respect to the distance function \eqref{pullbackdistance}.
\end{prop}
\begin{pf}
\begin{equation}
\begin{split}
&f^{-1}\left(\frac{\sum_{i=1}^{n}f\left(X_i^1\right)}{n}\right)=f^{-1}\left[\sum_{i=1}^{n}\frac{f\left(M\left(Z_{\neq i}(X_1^0,\dots,X_n^0)\right)\right)}{n}\right]=\\
&=f^{-1}\left[\sum_{i=1}^{n}\frac{\sum_{j=1,j\neq i}^{n-1}\frac{f\left(X_j^0\right)}{n-1}}{n}\right]=f^{-1}\left(\frac{\sum_{i=1}^{n}f\left(X_i^0\right)}{n}\right)
\end{split}
\end{equation}
Similarly we obtain the above equality for every iteration step, so
\begin{equation}
f^{-1}\left(\frac{\sum_{i=1}^{n}f\left(X_i^{l+1}\right)}{n}\right)=f^{-1}\left(\frac{\sum_{i=1}^{n}f\left(X_i^l\right)}{n}\right)=f^{-1}\left(\frac{\sum_{i=1}^{n}f\left(X_i^0\right)}{n}\right)\text{.}
\end{equation}
$\square$
\end{pf}

We turn our attention to more general functions than these pullback means. The following inequalities will turn out to be useful tools later.

\begin{lem}\label{kineqfunc}
\begin{equation}
\begin{split}
\left[(1-t)A^2+tB^2-\frac{k_2}{2}t(1-t)(A-B)^2\right]^{1/2}\leq\\
\leq \left[(1-t)A^2+tB^2-\frac{k_1}{2}t(1-t)(A-B)^2\right]^{1/2}
\end{split}
\end{equation}
for any $A,B\in P(r)$ and $t\in[0,1]$ if $k_1\leq k_2$.
\end{lem}
\begin{pf}
\begin{equation*}
\begin{split}
0\leq& \frac{k_2-k_1}{2}t(1-t)(A-B)^2\\
-\frac{k_2}{2}t(1-t)(A-B)^2\leq& -\frac{k_1}{2}t(1-t)(A-B)^2
\end{split}
\end{equation*}
\begin{equation*}
\begin{split}
(1-t)A^2+tB^2-\frac{k_2}{2}t(1-t)(A-B)^2&\leq\\
\leq &(1-t)A^2+tB^2-\frac{k_1}{2}t(1-t)(A-B)^2
\end{split}
\end{equation*}
Taking the square root of both sides and considering the fact that the the square root is operator monotone we get the inequality of the assertion.
$\square$
\end{pf}
Notice that for $k_1=0$ and $k_2=2$ we get the weighted arithmetic-square mean inequality
\begin{equation}\label{squarearithmetic}
(1-t)A+tB\leq \left[(1-t)A^2+tB^2\right]^{1/2}\text{.}
\end{equation}

We will prove an important inequality which will play a fundamental role in our further investigations. An important part of the proof of the convergence of the ALM- and BMP-process will rely on this inequality.

\begin{lem}\label{nonposcurv}
Let $k\in[0,2]$, $F:P(r)^2\mapsto P(r)$ and
\begin{equation}
F(A,B)\leq \left[(1-t)A^2+tB^2-\frac{k}{2}t(1-t)(A-B)^2\right]^{1/2}\text{.}
\end{equation}
Then with the distance function $d_E(A,B)^2=Tr\left\{(A-B)^*(A-B)\right\}$,
\begin{equation}
d_E\left(0,F(A,B)\right)^2\leq (1-t)d_E(0,A)^2+td_E(0,B)^2-\frac{k}{2}t(1-t)d_E(A,B)^2\text{.}
\end{equation}
\end{lem}
\begin{pf}
By substitution the assertion has the following form
\begin{equation}
Tr\left\{F(A,B)^2\right\}\leq Tr\left\{(1-t)A^2+tB^2-\frac{k}{2}t(1-t)(A-B)^2\right\}\text{.}
\end{equation}
This holds, since we have the following identity for hermitian positive definite $X\leq Y$
\begin{equation}
0\leq Tr\left\{Y^2-X^2\right\}=Tr\left\{(Y-X)(Y+X)\right\}\text{.}
\end{equation}
By choosing $X=F(A,B)-\left[(1-t)A^2+tB^2-\frac{k}{2}t(1-t)(A-B)^2\right]^{1/2}$ and $Y=F(A,B)+\left[(1-t)A^2+tB^2-\frac{k}{2}t(1-t)(A-B)^2\right]^{1/2}$ we get the assertion.
$\square$
\end{pf}

Notice that the above lemma is already true for every matrix mean $M(A,B)$ and their weighted $M_t(A,B)$ counterparts by Lemma~\ref{arithmeticmatrixmean} and \eqref{squarearithmetic}. By Lemma~\ref{kineqfunc} we also have a relatively wide family of functions which fulfills the conditions of the above lemma.

Now we are in position to prove Theorem~\ref{almthm}.
\begin{pf}[Theorem~\ref{almthm}]
The proof will be based on induction on the number of matrices $n$. We are going to measure the sum of the squared distances of the matrices $X_i^l$ from the zero matrix with respect to the distance function \eqref{frobeniusdistance} with
\begin{equation}
a_n^l=\sum_{i=1}^nd_E(0,X_i^l)^2=\sum_{i=1}^nTr\left\{(X_i^l)^2\right\}\text{.}
\end{equation}
We will also measure sum of the squared distances of the $X_i^l$ from one another. We will form this sum over all possible pairs of $X_i^l$ as
\begin{equation}
e_n^l=\sum_{1\leq i<j\leq n}d_E(X_i^l,X_j^l)^2=\sum_{1\leq i<j\leq n}Tr\left\{(X_i^l-X_j^l)^2\right\}\text{.}
\end{equation}

We are going to denote the common limit point of the sequences $X_i^l$ by $F_n(X_1^0,\ldots,X_n^0)$ for $n$. We will need the following lemmas which will be proved by induction as well on the number of matrices $n$, so we have to embed these lemmas into this proof of Theorem~\ref{almthm}. All three lemmas will be proved by assuming that they hold for $n$ matrices and also that the ALM-procedure converges to common limit for $n$ matrices. Making this assumption we show that the lemmas hold for $n+1$ and that the ALM procedure converges to common limit for $n+1$ as well. For the first step of the induction $(n=3)$ we will prove the lemmas directly. First we are going to treat the case of the first inequality \eqref{smaller}.
\begin{lem}[Monotone Iteration]
In the first case of inequality \eqref{smaller} we have
\begin{equation}
\left(\frac{\sum_{i=1}^{n}\left(X_i^{l+1}\right)^{-1}}{n}\right)^{-1}\geq \left(\frac{\sum_{i=1}^{n}\left(X_i^l\right)^{-1}}{n}\right)^{-1}\text{.}
\end{equation}
In the second case of inequality \eqref{arithmsmaller} we have
\begin{equation}
\frac{\sum_{i=1}^{n}X_i^{l+1}}{n}\leq \frac{\sum_{i=1}^{n}X_i^l}{n}\text{.}
\end{equation}
\end{lem}
\begin{pf}
We argue by induction on the number of matrices $n$ by making the assumption that the ALM-procedure converges for $n\geq 3$ to common limit point $F_n(X_1^0,\ldots,X_n^0)$, in other words $X_i^l\to F_n(X_1^0,\ldots,X_n^0)$ for $n$ and that the lemma holds for $n$. Consider the first case of inequality \eqref{smaller}. Then the inequality of the lemma for $n$ implies that
\begin{equation}\label{harmonicdecrease}
\left(\frac{\sum_{i=1}^{n}\left(X_i^l\right)^{-1}}{n}\right)^{-1}\geq \left(\frac{\sum_{i=1}^{n}\left(X_i^0\right)^{-1}}{n}\right)^{-1}
\end{equation}
and if we take the limit on the left hand side for $n$ we get the inequality
\begin{equation}\label{smaller2}
\left(\frac{\sum_{i=1}^{n}\left(F_n(X_1^0,\ldots,X_n^0)\right)^{-1}}{n}\right)^{-1}=F_n(X_1^0,\ldots,X_n^0)\geq \left(\frac{\sum_{i=1}^{n}\left(X_i^0\right)^{-1}}{n}\right)^{-1}\text{.}
\end{equation}
The above inequality also holds directly for $n=2$ by the assumption of inequality \eqref{smaller}, so this will also provide the first step in our induction.

Now we prove the lemma for $n+1$ if it is true for $n$. By \eqref{smaller2} we have
\begin{equation}
X_i^{l+1}=F_n\left(Z_{\neq i}(X_1^l,\ldots,X_{n+1}^l)\right)\geq \left(\frac{\sum_{j=1,j\neq i}^{n+1}\left(X_i^l\right)^{-1}}{n}\right)^{-1}\text{.}
\end{equation}
The $n+1$-variable harmonic mean is operator monotone in its variables, therefore if we take the $n+1$-variable harmonic mean of the above on the left and right hand side, we get
\begin{equation}
\left(\frac{\sum_{i=1}^{n+1}F_n\left(Z_{\neq i}(X_1^l,\ldots,X_{n+1}^l)\right)^{-1}}{n+1}\right)^{-1}\geq \left(\frac{\sum_{i=1}^{n+1}\left(H_i^l\right)^{-1}}{n+1}\right)^{-1}\text{,}
\end{equation}
where $H_i^l=\left(\frac{\sum_{j=1,j\neq i}^{n+1}\left(X_i^l\right)^{-1}}{n}\right)^{-1}$. Then by Proposition~\ref{invariantsquaremean} with $f(t)=t^{-1}$, the harmonic mean of the $n+1$ matrices is left invariant on the right hand side, so this is equivalent to
\begin{equation}
\left(\frac{\sum_{i=1}^{n+1}\left(X_i^{l+1}\right)^{-1}}{n+1}\right)^{-1}\geq \left(\frac{\sum_{i=1}^{n+1}\left(X_i^l\right)^{-1}}{n+1}\right)^{-1}\text{.}
\end{equation}


The second case given by inequality \eqref{arithmsmaller} is very similar to the proof of the first case. Instead of inequality \eqref{harmonicdecrease} we have
\begin{equation}
\frac{\sum_{i=1}^{n}X_i^l}{n}\leq \frac{\sum_{i=1}^{n}X_i^0}{n}
\end{equation}
and instead of \eqref{smaller2} we have
\begin{equation}
\frac{\sum_{i=1}^{n}F_n(X_1^0,\ldots,X_n^0)}{n}=F_n(X_1^0,\ldots,X_n^0)\leq \frac{\sum_{i=1}^{n}X_i^0}{n}\text{.}
\end{equation}
The rest of the argument is just the same, although we have the $n$-variable arithmetic mean replacing the $n$-variable harmonic mean, and the inequalities are reversed. The lemma is proved.
\end{pf}

\begin{lem}[Decreasing Distances]
We have
\begin{equation}\label{alman}
a_n^{l+1}\leq a_n^l-\frac{k}{8}z_ne_n^l
\end{equation}
in the case of \eqref{smaller}, or we have \eqref{alman} with $k=2$ in the case of \eqref{arithmsmaller}. In both cases $z_n=\frac{2}{n-1}$.
\end{lem}
\begin{pf}
We will prove this for the case \eqref{smaller}. The second case of \eqref{arithmsmaller} is just the same with $k=2$, we will only use that the right hand side of \eqref{smaller} holds, so we do not have to treat the second case \eqref{arithmsmaller} separately due to \eqref{squarearithmetic} with $t=1/2$. The first step is to show the above for $n=3$. By Lemma~\ref{nonposcurv} and that the right hand side of \eqref{smaller} is equivalent to the assumption of the lemma for $t=1/2$, we have
\begin{equation}
\begin{split}
d_E\left(0,F(X_i^l,X_j^l)\right)^2\leq \frac{d_E(0,X_i^l)^2+d(0,X_j^l)^2}{2}-\frac{k}{8}d_E(X_i^l,X_j^l)^2\\
d_E\left(0,X_s^{l+1}\right)^2\leq \frac{d_E(0,X_i^l)^2+d(0,X_j^l)^2}{2}-\frac{k}{8}d_E(X_i^l,X_j^l)^2\text{,}
\end{split}
\end{equation}
where $i,j,s\in\left\{1,2,3\right\}$ and $i\neq j\neq s, s\neq i$. There are 3 distinct inequalities of the above for $s=1,2,3$. By summing these inequalities for $s$ we get \eqref{alman} for $n=3$ and $z_3=1$.

Now suppose \eqref{alman} holds for $n$ and that $X_i^l$ converge to a common limit point for $n$ denoted again by $F_n(X_1^0,\ldots,X_n^0)$. Then we have
\begin{equation}
a_n^{l}\leq a_n^0-\frac{k}{8}z_ne_n^0
\end{equation}
and by taking the limit on the left hand side we get
\begin{equation}\label{inequdecrease}
\begin{split}
\lim_{l\to\infty}a_n^{l}=nd_E\left(0,F_n(X_1^0,\ldots,X_n^0)\right)^2\leq a_n^0-\frac{k}{8}z_ne_n^0\\
d_E\left(0,F_n(X_1^0,\ldots,X_n^0)\right)^2\leq \frac{a_n^0}{n}-\frac{k}{8}\frac{z_n}{n}e_n^0\text{.}
\end{split}
\end{equation}
Then set up the ALM-procedure on $X_i^0\in P(r); i=1,2\ldots,n+1$ with $M_n(X_1,\ldots,X_n):=F_n(X_1,\ldots,X_n)$. Inequality \eqref{inequdecrease} can be applied in any of the iteration steps, so we get
\begin{equation}
\begin{split}
&d_E\left(0,F_n\left(Z_{\neq i}\left(X_1^l,\ldots,X_{n+1}^l\right)\right)\right)^2\leq \\
&\leq\frac{\sum_{j=1,j\neq i}^{n+1}d_E(0,X_j^l)^2}{n}-\frac{k}{8}\frac{z_n}{n}\sum_{1\leq j<s\leq n+1,j\neq i,s\neq i}d_E(X_j^l,X_s^l)^2\text{.}
\end{split}
\end{equation}
If we sum these inequalities for $i$ we arrive at the following
\begin{equation}
\begin{split}
&\sum_{i=1}^{n+1}d_E\left(0,F_n\left(Z_{\neq i}\left(X_1^l,\ldots,X_{n+1}^l\right)\right)\right)^2\leq \\
&\leq\sum_{i=1}^{n+1}\frac{\sum_{j=1,j\neq i}^{n+1}d_E(0,X_j^l)^2}{n}-\frac{k}{8}\frac{z_n}{n}\sum_{i=1}^{n+1}\sum_{1\leq j<s\leq n+1,j\neq i,s\neq i}d_E(X_j^l,X_s^l)^2\text{.}
\end{split}
\end{equation}
The left hand side of the above is just $a_{n+1}^{l+1}$. The first term on the right hand side is easily written as
\begin{equation}
\sum_{i=1}^{n+1}\frac{\sum_{j=1,j\neq i}^{n+1}d_E(0,X_j^l)^2}{n}=\sum_{i=1}^{n+1}d_E(0,X_i^l)^2=a_{n+1}^l\text{.}
\end{equation}
We have to carefully analyze the second term
\begin{equation}\label{subgraphsum}
\sum_{i=1}^{n+1}\sum_{1\leq j<s\leq n+1,j\neq i,s\neq i}d_E(X_j^l,X_s^l)^2\text{.}
\end{equation}
Consider the complete graph $K_{n+1}$ on $n+1$ vertices labelled from $1$ to $n+1$. In this way we have a natural bijective mapping between the matrices $X_i^l$ and the vertices of $K_{n+1}$. Then for every squared distance $d_E(X_j^l,X_s^l)^2$ we have a corresponding edge in $K_{n+1}$ of the form $(j,s)$. Then the sum $\sum_{1\leq j<s\leq n+1,j\neq i,s\neq i}d_E(X_j^l,X_s^l)^2$ is just the sum of the squared distances corresponding to the edges of the complete graph $K_n$ given on the vertices $\left\{1,\ldots,i-1,i+1,\ldots,n+1\right\}$. This is almost $\sum_{1\leq j<s\leq n+1}d_E(X_j^l,X_s^l)^2$, but we leave out from the sum every squared distance corresponding to an edge that has the vertex $i$ as an ending vertex. So actually \eqref{subgraphsum} almost equals to
\begin{equation}
\sum_{i=1}^{n+1}\sum_{1\leq j<s\leq n+1}d_E(X_j^l,X_s^l)^2=(n+1)e_{n+1}^l\text{,}
\end{equation}
but in the sum \eqref{subgraphsum} every vertex has been left out once, so every squared distance corresponding to an edge has been left out twice, hence
\begin{equation}\label{graphsum}
\sum_{i=1}^{n+1}\sum_{1\leq j<s\leq n+1,j\neq i,s\neq i}d_E(X_j^l,X_s^l)^2=(n-1)e_{n+1}^l\text{.}
\end{equation}
This shows us that $z_{n+1}=\frac{n-1}{n}z_n$ and also $z_3=1$, so in other words by solving the recursion we get
\begin{equation}
z_n=\frac{2}{n-1}\text{.}
\end{equation}
This concludes the lemma for every $n$.
\end{pf}
\begin{lem}[Boundedness]
The matrix sequences $X_i^l$ are bounded for all $n$.
\end{lem}
\begin{pf}
We have the trivial lower bound $X_i^l\geq 0$, since by assumption $F(A,B)\geq 0$, so we have $X_i^l\geq 0$ for $n=3$. Now assume again that the ALM-procedure converges to common limit point denoted by $F_n(X_1^0,\ldots,X_n^0)$ for $n$ and $F_n(X_1^0,\ldots,X_n^0)\geq 0$. Then trivially for $n+1$ the sequences $X_i^l\geq 0$ since $F_n(X_1^0,\ldots,X_n^0)\geq 0$. This also shows that if the sequences converge for $n+1$ to a common limit $F_{n+1}(X_1^0,\ldots,X_n^0)$, then this limit is also bounded from below, so $F_{n+1}(X_1^0,\ldots,X_n^0)\geq 0$.

Now we provide a suitable upper bound as well. By the previous assertion Lemma Decreasing Distances we have \eqref{alman} for $n\geq 3$. In particularly for $n=3$ it holds providing the first step, while for $n>3$ we need the inductional hypothesis that the ALM-procedure converges to a common limit point $F_{n-1}(X_1^0,\ldots,X_{n-1}^0)$ for $n-1$. The rest of the argument is just the same for all $n\geq 3$. So by \eqref{alman} we have $a_{n}^{l+1}\leq a_{n}^{l}$ which means that the sequence is monotone decreasing. So we have the bound
\begin{equation}
a_{n}^{l}\leq a_{n}^{0}=b\text{.}
\end{equation}
By this above we get for arbitrary $i$ that
\begin{equation}
d\left(0,X_i^l\right)^2=a_{n}^{l}-\sum_{j=1,j\neq i}^{n}d\left(0,X_j^l\right)^2\leq a_{n}^{l}\leq b\text{.}
\end{equation}
This means that the norm $\left\|X_i^l\right\|$ is bounded from above by $b$ since
\begin{equation}
\left\|X_i^l\right\|^2=Tr\left\{\left(X_i^l\right)^2\right\}=d\left(0,X_i^l\right)^2\text{.}
\end{equation}
This concludes the proof of the lemma.
\end{pf}

Now we move on to the final step of the induction. We prove that for $n=3$ the ALM-procedure converges and that if it converges for $n$ then it converges for $n+1$ in both cases of inequalities \eqref{smaller} and \eqref{arithmsmaller}. This last step will be based on the three lemmas: Lemma Monotone Iteration, Lemma Decreasing Distances and Lemma Boundedness. It is not necessary to prove separately the $n=3$ case since these three lemmas hold for $n=3$ and the argument will be the same as for general $n+1$ requiring the inductional hypothesis, the convergence of the procedure to a common limit point for $n$.

So by Lemma Decreasing Distances we have
\begin{equation}
a_n^{l+1}\leq a_n^l-\frac{k}{8}z_ne_n^l\text{,}
\end{equation}
in other words $a_n^l$ is a decreasing nonnegative sequence in $l$, therefore convergent. Since $z_n>0$ and has fixed value for each $n$ by Lemma Decreasing Distances, this means that $e_n^l\to 0$ as $l\to\infty$, so the matrices $X_i^l$ are approaching one another. By Lemma Boundedness we have that these sequences are bounded, hence they have convergent subsequences. But since $e_n^l\to 0$ these subsequences are converging to a common limit point. Let $X_i^{s_l}$ denote a subsequence converging to say $A$ and $X_i^{r_l}$ another subsequence converging to $B$. Without loss of generality we can take $s_l>r_l$. By Lemma Monotone Iteration for the case of inequality \eqref{smaller} we have
\begin{equation}
\left(\frac{\sum_{i=1}^{n}\left(X_i^{s_l}\right)^{-1}}{n}\right)^{-1}\geq \left(\frac{\sum_{i=1}^{n}\left(X_i^{r_l}\right)^{-1}}{n}\right)^{-1}\text{.}
\end{equation}
Now choose a subsequence of subsequences $s_j<r_j$ so then again by the lemma
\begin{equation}
\left(\frac{\sum_{i=1}^{n}\left(X_i^{s_j}\right)^{-1}}{n}\right)^{-1}\leq \left(\frac{\sum_{i=1}^{n}\left(X_i^{r_j}\right)^{-1}}{n}\right)^{-1}\text{.}
\end{equation}
Taking the limits we have $A\geq B$ and $A\leq B$ so $A=B$. In the second case \eqref{arithmsmaller} we have the same argument but using the $n$-variable arithmetic mean instead of the $n$-variable harmonic above.

Now this argument shows the convergence of the ALM-procedure to a common limit point directly for $n=3$ and inductively for $n$ assuming convergence to common limit for $n-1$.
$\square$
\end{pf}

We are going to study some properties of this limit point later, jointly with the case of the BMP-mean, after showing that the BMP procedure converges. In the next section we will show a similar theorem to Theorem~\ref{almthm} for the BMP-procedure.

\section{Bini-Meini-Poloni procedure for every matrix mean}
In this section we will treat the case of the Bini-Meini-Poloni procedure. We may do that for matrix means since we have defined a weighted mean $M_t(A,B)$ corresponding to any symmetric matrix mean $M(A,B)$. The outline of the proof of the convergence of the BMP-procedure will roughly follow the one of the ALM-procedure, although some lemmas will be formulated differently.

Firstly let us recall the Bini-Meini-Poloni procedure \cite{bini}:

\begin{defn}\label{bmpit}[BMP iteration]
Let $X=(X_1^0,\dots,X_n^0)$ where $X_i^0 \in P(r)$ and define the mapping $M(X_1,\dots,X_{n})$ inductively as follows. If $n=2$ assume that $M_t(X_1,X_2)$ is already given. For general $n>2$ assume that $M(X_1,\dots,X_{n-1})$ is already defined. Then using $M(X_1,\dots,X_{n-1})$, set up the iteration
\begin{equation}
X_i^{l+1}=M_{\frac{n-1}{n}}\left(X_i^l,M\left(Z_{\neq i}\left(X_1^l,\dots,X_n^l\right)\right)\right)\text{,}
\end{equation}
where $Z_{\neq i}(X_1^l,\dots,X_n^l)=X_1^l,\dots,X_{i-1}^l, X_{i+1}^l,\dots,X_n^l$. If the sequences $X_i^{l}$ converge to a common limit point for every $i$, then define
\begin{equation}
\lim_{l\to\infty}X_i^{l}=M(X_1^0,\dots,X_n^0)\text{.}
\end{equation}
\end{defn}

\begin{thm}\label{bmpthm}
Let $F:[0,1]\times P(r)^2\mapsto P(r)$ and suppose that $F_t(A,B)$ fulfills one of the inequalities below:
\begin{equation}\label{smallerbmp}
\begin{split}
&\left[(1-t)A^{-1}+tB^{-1}\right]^{-1}\leq F_t(A,B)\leq\\
&\leq\left[(1-t)A^2+tB^2-\frac{k}{2}t(1-t)(A-B)^2\right]^{1/2}
\end{split}
\end{equation}
for a $k\in(0,2]$ and every $t\in [0,1]$, or
\begin{equation}\label{arithmsmallerbmp}
F_t(A,B)\leq (1-t)A+tB\text{,}
\end{equation}
for every $t\in [0,1]$. Then in Definition~\ref{bmpit} starting with $M_t(A,B):=F_t(A,B)$, $M(X_1,\dots, X_{n})$ exists for all $n$, in other words the sequences converge to a common limit point for all $n$.
\end{thm}

Before we turn to the proof of the above theorem, we again consider some lemmas which will be similar to the ALM case. Let us recall again the metric space $P(r)$ with the distance function \eqref{pullbackdistance}. We already know that the minimum of
\begin{equation}
b(x)=\sum_{i=1}^nd_f(x,x_i)^2
\end{equation}
is attained at a unique point in $P(r)$ denoted by $\hat{x}$ and we also know that
\begin{equation}
\hat{x}=f^{-1}\left(\frac{\sum_{i=1}^nf(x_i)}{n}\right)\text{.}
\end{equation}

We will need a similar theorem to Proposition~\ref{invariantsquaremean}.
\begin{prop}\label{invariantsquaremeanbmp}
Let $X_i^0\in P(r)$ for $i=1,\ldots,n$. Then the BMP-procedure (Definition~\ref{bmpit}) set up on the matrices $X_1^0,\dots,X_n^0$ with the weighted mean function $M_t(A,B):=f^{-1}\left((1-t)f(A)+tf(B)\right)$ and the $n-1$ variable function $M(x_1,\dots,x_{n-1}):=f^{-1}\left(\frac{\sum_{i=1}^{n-1}f(x_i)}{n-1}\right)$ leaves the Riemann centroid of the points $X_1^0,\dots,X_n^0$ invariant with respect to the distance function \eqref{pullbackdistance}.
\end{prop}
\begin{pf}
\begin{equation}
\begin{split}
&f^{-1}\left(\frac{\sum_{i=1}^{n}f\left(X_i^1\right)}{n}\right)=f^{-1}\left[\sum_{i=1}^{n}\frac{f\left(M_{\frac{n-1}{n}}\left(X_i^0,M\left(Z_{\neq i}(X_1^0,\dots,X_n^0)\right)\right)\right)}{n}\right]=\\
&=f^{-1}\left[\sum_{i=1}^n\frac{\frac{f(X_i^0)}{n}+\frac{n-1}{n}\sum_{j=1,j\neq i}^n\frac{f(X_j^0)}{n-1}}{n}\right]=f^{-1}\left(\frac{\sum_{i=1}^{n}f\left(X_i^0\right)}{n}\right)
\end{split}
\end{equation}
Similarly we obtain the above equality for every iteration step, so
\begin{equation}
f^{-1}\left(\frac{\sum_{i=1}^{n}f\left(X_i^{l+1}\right)}{n}\right)=f^{-1}\left(\frac{\sum_{i=1}^{n}f\left(X_i^l\right)}{n}\right)=f^{-1}\left(\frac{\sum_{i=1}^{n}f\left(X_i^0\right)}{n}\right)\text{.}
\end{equation}
$\square$
\end{pf}

\begin{pf}[Theorem~\ref{bmpthm}]
The proof again will be based on induction on the number of matrices $n$. We will use the same notations to denote the sum of the squared distances of the matrices $X_i^l$ from the zero matrix with respect to the distance function \eqref{frobeniusdistance}, so
\begin{gather}
a_n^l=\sum_{i=1}^nd_E(0,X_i^l)^2=\sum_{i=1}^nTr\left\{(X_i^l)^2\right\}\\
e_n^l=\sum_{1\leq i<j\leq n}d_E(X_i^l,X_j^l)^2=\sum_{1\leq i<j\leq n}Tr\left\{(X_i^l-X_j^l)^2\right\}\text{.}
\end{gather}

We will denote by $F(X_1^0,\ldots,X_n^0)$ the common limit point of the sequences $X_i^l$ for $n$. The proof will rely on similar three lemmas to the ones in the proof of the ALM-procedure. First we are going to treat the case of the first inequality \eqref{smallerbmp}.
\begin{lem}[Monotone Iteration]
In the first case of inequality \eqref{smallerbmp} we have
\begin{equation}
\left(\frac{\sum_{i=1}^{n}\left(X_i^{l+1}\right)^{-1}}{n}\right)^{-1}\geq \left(\frac{\sum_{i=1}^{n}\left(X_i^l\right)^{-1}}{n}\right)^{-1}\text{.}
\end{equation}
In the second case of inequality \eqref{arithmsmallerbmp} we have
\begin{equation}
\frac{\sum_{i=1}^{n}X_i^{l+1}}{n}\leq \frac{\sum_{i=1}^{n}X_i^l}{n}\text{.}
\end{equation}
\end{lem}
\begin{pf}
The proof uses similar ideas to the case of Lemma Monotone Iteration for the ALM-process. We again argue by induction on the number of matrices $n$. Consider the first case of inequality \eqref{smallerbmp}. Suppose that the BMP-procedure converges for $n\geq 3$ to common limit point $F(X_1^0,\ldots,X_n^0)$ in other words $X_i^l\to F(X_1^0,\ldots,X_n^0)$ for $n$. Also the inequality of the lemma for $n$ implies that
\begin{equation}\label{harmonicdecreasebmp}
\left(\frac{\sum_{i=1}^{n}\left(X_i^l\right)^{-1}}{n}\right)^{-1}\geq \left(\frac{\sum_{i=1}^{n}\left(X_i^0\right)^{-1}}{n}\right)^{-1}
\end{equation}
and if we take the limit on the left hand side for $n$ we get the inequality
\begin{equation}\label{smaller2bmp}
\left(\frac{\sum_{i=1}^{n}\left(F(X_1^0,\ldots,X_n^0)\right)^{-1}}{n}\right)^{-1}=F(X_1^0,\ldots,X_n^0)\geq \left(\frac{\sum_{i=1}^{n}\left(X_i^0\right)^{-1}}{n}\right)^{-1}\text{.}
\end{equation}
The above inequality also holds for $n=2$ by the assumption of inequality \eqref{smallerbmp}, so this provides the first step for $n=3$ in our induction.

Now we prove the lemma for $n+1$ if it is true for $n$. Similarly to the case of Lemma Monotone Iteration in the ALM process, we make use of the operator monotonicity of the $n+1$-variable harmonic mean, and use \eqref{smaller2bmp}. Then we use Proposition~\ref{invariantsquaremeanbmp} with $f(t)=t^{-1}$ and the same argument as in the ALM case, performed using instead one BMP iteration step, yields
\begin{equation}
\left(\frac{\sum_{i=1}^{n+1}\left(X_i^{l+1}\right)^{-1}}{n+1}\right)^{-1}\geq \left(\frac{\sum_{i=1}^{n+1}\left(X_i^l\right)^{-1}}{n+1}\right)^{-1}\text{.}
\end{equation}


The second case given by inequality \eqref{arithmsmallerbmp} again can be treated similarly to the ALM-case.
\end{pf}

\begin{lem}[Decreasing Distances]
We have
\begin{equation}\label{bmpan}
a_n^{l+1}\leq a_n^l-\frac{k}{8}z_ne_n^l\text{,}
\end{equation}
in the case of \eqref{smallerbmp}, or we have \eqref{bmpan} with $k=2$ in the case of \eqref{arithmsmallerbmp}. In both cases $z_n=\frac{4}{(n-1)n}$.
\end{lem}
\begin{pf}
The situation is similar again to the ALM case. We will prove this for the case \eqref{smallerbmp}, the second case of \eqref{arithmsmallerbmp} is just the same with $k=2$, since we will only use that the right hand side of \eqref{smaller} holds, so we do not have to treat the second case \eqref{arithmsmallerbmp} separately due to \eqref{squarearithmetic}. The first step is to show the above for $n=3$. By Lemma~\ref{nonposcurv} and that the right hand side of \eqref{smallerbmp} is equivalent to the assumption of the lemma, we get
\begin{equation}
\begin{split}
d_E\left(0,F(X_i^l,X_j^l)\right)^2\leq& \frac{d_E(0,X_i^l)^2+d(0,X_j^l)^2}{2}-\frac{k}{8}d_E(X_i^l,X_j^l)^2\\
d_E\left(0,F_{2/3}\left(X_s^l,F(X_i^l,X_j^l)\right)\right)^2\leq& \frac{1}{3}d_E(0,X_s^l)^2+\frac{2}{3}d(0,F(X_i^l,X_j^l))^2-\\
&-\frac{k}{2}\frac{1}{3}\frac{2}{3}d_E(X_s^l,F(X_i^l,X_j^l))^2\\
d_E\left(0,F_{2/3}\left(X_s^l,F(X_i^l,X_j^l)\right)\right)^2\leq& \frac{1}{3}d_E(0,X_s^l)^2+\frac{2}{3}d(0,F(X_i^l,X_j^l))^2-\\
d_E\left(0,F_{2/3}\left(X_s^l,F(X_i^l,X_j^l)\right)\right)^2\leq& \frac{d_E(0,X_s^l)^2+d_E(0,X_i^l)^2+d(0,X_j^l)^2}{3}-\\
&-\frac{k}{8}\frac{2}{3}d_E(X_i^l,X_j^l)^2
\end{split}
\end{equation}
in other words the last inequality is equivalent to
\begin{equation}
d_E\left(0,X_s^{l+1}\right)^2\leq \frac{d_E(0,X_s^l)^2+d_E(0,X_i^l)^2+d(0,X_j^l)^2}{3}-\frac{k}{8}\frac{2}{3}d_E(X_i^l,X_j^l)^2\text{,}
\end{equation}
where $i,j,s\in\left\{1,2,3\right\}$ and $i\neq j\neq s, s\neq i$. There are 3 distinct inequalities of the above for $s=1,2,3$. By summing these inequalities for $s$ we get \eqref{bmpan} for $n=3$ and $z_3=\frac{2}{3}$.

Now suppose \eqref{bmpan} holds for $n$ and that $X_i^l$ converge to a common limit point for $n$. Then we have
\begin{equation}
a_n^{l}\leq a_n^0-\frac{k}{8}z_ne_n^0
\end{equation}
and by taking the limit on the left hand side we get
\begin{equation}\label{inequbmp}
\begin{split}
\lim_{l\to\infty}a_n^{l}=nd_E\left(0,F(X_1^0,\ldots,X_n^0)\right)^2\leq a_n^0-\frac{k}{8}z_ne_n^0\\
d_E\left(0,F(X_1^0,\ldots,X_n^0)\right)^2\leq \frac{a_n^0}{n}-\frac{k}{8}\frac{z_n}{n}e_n^0\text{.}
\end{split}
\end{equation}
Then set up the BMP-procedure on $X_i^0\in P(r); i=1,2\ldots,n+1$ with $M_n(X_1,\ldots,X_n):=F(X_1,\ldots,X_n)$. Inequality \eqref{inequbmp} can be applied in any of the iteration steps, so we get
\begin{equation}\label{inequbmp2}
\begin{split}
&d_E\left(0,F\left(Z_{\neq i}\left(X_1^l,\ldots,X_{n+1}^l\right)\right)\right)^2\leq \\
&\leq\frac{\sum_{j=1,j\neq i}^{n+1}d_E(0,X_j^l)^2}{n}-\frac{k}{8}\frac{z_n}{n}\sum_{1\leq j<s\leq n+1,j\neq i,s\neq i}d_E(X_j^l,X_s^l)^2\text{.}
\end{split}
\end{equation}
Then we have to compute $X_i^{l+1}=F_{\frac{n}{n+1}}\left(X_i^l,F\left(Z_{\neq i}\left(X_1^l,\ldots,X_{n+1}^l\right)\right)\right)$ and bound its squared distance from the zero matrix
\begin{equation}
\begin{split}
d_E\left(0,F_{\frac{n}{n+1}}\left(X_i^l,F\left(Z_{\neq i}\left(X_1^l,\ldots,X_{n+1}^l\right)\right)\right)\right)^2\leq \frac{1}{n+1}d_E\left(0,X_i^l\right)^2+\\
+\frac{n}{n+1}d_E\left(0,F\left(Z_{\neq i}\left(X_1^l,\ldots,X_{n+1}^l\right)\right)\right)^2-\\
-\frac{k}{8}\frac{1}{n+1}\frac{n}{n+1}d_E\left(X_i^l,F_{\frac{n}{n+1}}\left(X_i^l,F\left(Z_{\neq i}\left(X_1^l,\ldots,X_{n+1}^l\right)\right)\right)\right)^2\text{.}
\end{split}
\end{equation}
We drop the last term, as it seems that it is hard to estimate it from below, and substitute in using inequality \eqref{inequbmp2}, we get
\begin{equation}
\begin{split}
&d_E\left(0,X_i^{l+1}\right)^2\leq \frac{1}{n+1}d_E\left(0,X_i^l\right)^2+\frac{n}{n+1}d_E\left(0,F\left(Z_{\neq i}\left(X_1^l,\ldots,X_{n+1}^l\right)\right)\right)^2\\
&\leq \frac{1}{n+1}d_E\left(0,X_i^l\right)^2+\\
&+\frac{n}{n+1}\left[\frac{\sum_{j=1,j\neq i}^{n}d_E(0,X_j^l)^2}{n}-\frac{k}{8}\frac{z_n}{n}\sum_{1\leq j<s\leq n+1,j\neq i,s\neq i}d_E(X_j^l,X_s^l)^2\right]\\
&\leq \frac{\sum_{j=1,j\neq i}^{n+1}d_E(0,X_j^l)^2}{n+1}-\frac{k}{8}\frac{z_n}{n+1}\sum_{1\leq j<s\leq n+1,j\neq i,s\neq i}d_E(X_j^l,X_s^l)^2\text{.}
\end{split}
\end{equation}
If we sum these inequalities for $i$ we arrive at the following
\begin{equation}
\begin{split}
&\sum_{i=1}^{n+1}d_E\left(0,X_i^{l+1}\right)^2\leq \\
&\leq\sum_{i=1}^{n+1}\frac{\sum_{j=1}^{n+1}d_E(0,X_j^l)^2}{n+1}-\frac{k}{8}\frac{z_n}{n+1}\sum_{i=1}^{n+1}\sum_{1\leq j<s\leq n+1,j\neq i,s\neq i}d_E(X_j^l,X_s^l)^2\text{,}
\end{split}
\end{equation}
which is equivalent to
\begin{equation}
\begin{split}
&\sum_{i=1}^{n+1}d_E\left(0,X_i^{l+1}\right)^2\leq \\
&\leq\sum_{i=1}^{n+1}d_E(0,X_i^l)^2-\frac{k}{8}\frac{z_n}{n+1}\sum_{i=1}^{n+1}\sum_{1\leq j<s\leq n+1,j\neq i,s\neq i}d_E(X_j^l,X_s^l)^2\text{.}
\end{split}
\end{equation}
The left hand side of the above is just $a_{n+1}^{l+1}$ and the first term on the right hand side is $a_{n+1}^l$. By the proof of the convergence of the ALM-process the second term
\begin{equation}
\sum_{i=1}^{n+1}\sum_{1\leq j<s\leq n+1,j\neq i,s\neq i}d_E(X_j^l,X_s^l)^2=(n-1)e_{n+1}^l\text{.}
\end{equation}
Thus $z_{n+1}=\frac{n-1}{n+1}z_n$ and also $z_3=\frac{2}{3}$, by solving the recursion we get
\begin{equation}
z_n=\frac{4}{(n-1)n}\text{.}
\end{equation}
This concludes the lemma for every $n$.
\end{pf}
\begin{lem}[Boundedness]
The matrix sequences $X_i^l$ are bounded for all $n$.
\end{lem}
\begin{pf}
We have the trivial lower bound $X_i^l\geq 0$, since by assumption $F_t(A,B)\geq 0$, so we have $X_i^l\geq 0$ for $n=3$. Now similarly to the case of the ALM-process we assume again that the BMP-procedure converges to common limit point denoted by $F(X_1^0,\ldots,X_n^0)$ for $n$ and $F(X_1^0,\ldots,X_n^0)\geq 0$. Then again if the sequences converge for $n+1$ to a common limit $F(X_1^0,\ldots,X_{n+1}^0)$, then this limit $F(X_1^0,\ldots,X_{n+1}^0)\geq 0$.

We provide the suitable upper bound similarly to the ALM case. We again have $a_{n+1}^{l}\leq a_{n+1}^{0}=b$ so we similarly get $\left\|X_i^l\right\|\leq b$ for $n+1$ if the procedure converges to common limit for $n$. This finishes the proof of the lemma.
\end{pf}

Now the last step of the proof is exactly the same as in the case of the ALM-procedure, we just have to use the three lemmas: Lemma Monotone Iteration, Decreasing Distances and Boundedness adapted for the case of the BMP iteration.
$\square$
\end{pf}

Notice that for matrix means we have the weighted mean procedure $M_t(A,B)$ introduced here. By Proposition~\ref{weightedmeanproperties} we have that every such mean is smaller than the weighted arithmetic mean and larger than the weighted harmonic mean. So as a consequence of Theorem~\ref{bmpthm} we get that the BMP-procedure converges for every symmetric matrix mean if we identify their weighted counterparts with our weighted mean $M_t(A,B)$.

\section{Properties of the ALM and BMP mean}
We will show that the limit point of the ALM and BMP processes, denoted by $M_{ALM}(X_1,\ldots,X_n)$ and $M_{BMP}(X_1,\ldots,X_n)$ respectively, as extensions of symmetric matrix means, fulfill the following properties.
\begin{thm}\label{properties}
If $M(A,B)$ is a symmetric matrix mean, then the $M:=M_{ALM}(X_1,\ldots,X_n)$ and $M:=M_{BMP}(X_1,\ldots,X_n)$ extensions fulfill the following properties
\begin{enumerate}
\renewcommand{\labelenumi}{(\Roman{enumi})}
\item $M(X,\dots, X)=X$ for every $X\in P(r)$,
\item $M(X_1,\ldots,X_n)$ is invariant under the permutation of its variables,
\item $\min(X_1,\dots, X_n)\leq M(X_1,\dots, X_n)\leq \max(X_1,\dots, X_n)$ if $\min$ and $\max$ exist with respect to the positive definite order,
\item If $X_i\leq X'_i$, then $M(X_1,\dots, X_n)\leq M(X'_1,\dots, X'_n)$,
\item $M(X_1,\ldots,X_n)$ is continuous,
\item $M(CX_1C^{*},\dots, CX_nC^{*})=CM(X_1,\dots, X_n)C^{*}\text{ for all invertible C}$.
\end{enumerate}
\end{thm}
\begin{pf}
The proof of each property will be based on induction. Each of them trivially holds for $n=2$ by properties of matrix means discussed in the second section. So it remains to prove them for $n+1$ assuming that they hold for $n$.

Property (I) and (II) trivially holds for $n+1$ if it holds for $n$. We prove property (IV). Let $X_{1}^0,\ldots, X_{n+1}^0 \in P(r)$ and $X_{i}^0\leq (X'_{i})^0\in P(r)$. If we iterate by the ALM process, it is easy to see that the order $X_{i}^0\leq (X'_{i})^0$ is preserved due to the inductional hypothesis on property (IV), so $X_{i}^l\leq (X'_{i})^l$. Taking the limits $l\to\infty$ we get the assertion. In case of the BMP-process the argument is similar but we have to use also that $M_t(A,B)\leq M_t(A',B')$ if $A\leq A'$ and $B\leq B'$.

Property (III) is an easy consequence of property (I) and (IV), if minimum and maximum exist. Setting up the same iteration on the new $n$-tuple formed by the minimal element we get the inequality on the left in property (III), similarly we can obtain the inequality on the right as well.

To prove property (VI) let $(X'_i)^{0}=CX_{i}^{0}C^{*}$ and set up the ALM or BMP process on $(X'_1)^{0},\ldots, (X'_n)^{0}$ as on $X_1^{0},\ldots, X_n^{0}$. Property (VI) implies in the case of ALM
\begin{equation}
\begin{split}
CX_{i}^{l+1}C^{*}=CM\left(Z_{\neq i}\left(X_{1}^{l},\ldots,X_{n+1}^{l}\right)\right)C^{*}=\\
=M\left(Z_{\neq i}\left(CX_{1}^{l}C^{*},\ldots,CX_{n+1}^{l}C^{*}\right)\right)\text{.}
\end{split}
\end{equation}
In the case of the BMP process we have similarly
\begin{equation}
\begin{split}
CX_{i}^{l+1}C^{*}=CM_{\frac{n}{n+1}}\left(X_i^{l},M\left(Z_{\neq i}\left(X_{1}^{l},\ldots,X_{n+1}^{l}\right)\right)\right)C^{*}=\\
=M_{\frac{n}{n+1}}\left(CX_i^{l}C^{*},M\left(Z_{\neq i}\left(CX_{1}^{l}C^{*},\ldots,CX_{n+1}^{l}C^{*}\right)\right)\right)\text{.}
\end{split}
\end{equation}
Applying the above recursively in every iteration step we get
\begin{equation}
CX_{i}^{l}C^{*}=(X'_i)^{l}\text{.}
\end{equation}
Taking the limit $l\to \infty$ the assertion follows.

Property (V) is a consequence of properties (IV) and (VI) by Lemma~\ref{monotonecont}.
$\square$
\end{pf}
We also have that the ALM and BMP procedures preserve the ordering of functions. So we have for the ALM process the following
\begin{prop}
If $M(A,B)\leq N(A,B)$ are functions satisfying the properties of $F(A,B)$ in Theorem~\ref{almthm}, then the same ordering is true for the ALM limit points $M(X_1,\ldots, X_n)$ and $N(X_1,\ldots, X_n)$.
\end{prop}
\begin{pf}
Again we argue by induction. The inequality $M(X_1,\ldots, X_{n})\leq N(X_1,\ldots, X_{n})$ holds for $n=2$ by assumption. Let us denote the matrices in the ALM iteration steps performed with $M(X_1,\ldots, X_{n-1})$ and $N(X_1,\ldots, X_{n-1})$ on $X_1^0,\ldots,X_{n}^0\in P(r)$ by $X_i^l$ and $(X_i')^l$ respectively. Now again we have $M(X_1,\ldots, X_{n-1})\leq N(X_1,\ldots, X_{n-1})$ by the inductional hypothesis so we have $X_i^l\leq (X_i')^l$. Taking the limits we get the assertion.
$\square$
\end{pf}
A similar, although a bit different assertion holds for the BMP process.
\begin{prop}
If $M_t(A,B)\leq N_t(A,B)$ are functions satisfying the properties of $F_t(A,B)$ in Theorem~\ref{bmpthm}, then the same ordering is true for the BMP limit points $M(X_1,\ldots, X_n)$ and $N(X_1,\ldots, X_n)$.
\end{prop}
\begin{pf}
We have an inductional argument similarly to the preceding case of the ALM process. The inequality $M(A,B)\leq N(A,B)$ holds for by assumption since $M(A,B)=M_{1/2}(A,B)$ and $N(A,B)=N_{1/2}(A,B)$. Let us denote the matrices in the BMP iteration steps performed with $M(X_1,\ldots, X_{n-1})$ and $N(X_1,\ldots, X_{n-1})$ on $X_1^0,\ldots,X_{n}^0\in P(r)$ by $X_i^l$ and $(X_i')^l$ respectively. Now again we have $M(X_1,\ldots, X_{n-1})\leq N(X_1,\ldots, X_{n-1})$ by the inductional hypothesis and also $M_t(A,B)\leq N_t(A,B)$ so we have $X_i^l\leq (X_i')^l$. Taking the limits we get the assertion.
$\square$
\end{pf}

In the next section we will consider some convergence rate properties fulfilled by the BMP process.

\section{Convergence rate of the BMP process}
A well known property of the BMP process is its cubic convergence rate in a small neighborhood of its limit point. This is an advantage over the ALM process which is known to converge linearly. Or to be more specific these convergence rates are only precisely known in the case of the geometric mean \cite{ando, bini}. In this section we will show that the BMP process generally converges cubically for every possible matrix mean in a small neighborhood of the limit point of the process. The proof will be similar to the one presented in \cite{bini}. In order to be able to use such an argument we have to obtain a series expansion for the weighted mean $M_t(A,B)$ in the neighborhood of the identity matrix $I$.

We are going to use the big O notation as it is used in \cite{bini}. This means that we have $X=Y+O(\epsilon^k)$ if and only if there exist constants $\epsilon_0<1$ and $\theta$ such that for each $0<\epsilon<\epsilon_0$ we have $\left\|X-Y\right\|\leq \theta\epsilon^k$.
\begin{prop}\label{expansion}
Let $M(A,B)$ be a symmetric matrix mean and $f(t)$ be its corresponding normalized operator monotone function. Let $f(t)$ have a series expansion around $I$ as
\begin{equation}\label{fexpansion}
f(X)=I+\frac{X-I}{2}+\sum_{k=2}^{\infty}b_k(X-I)^{k}\text{.}
\end{equation}
Then we have a series expansion for $M_t(I,X)=f_{t}(X)$ whenever $\left\|X-I\right\|\leq\epsilon<1$ in the form
\begin{equation}
f_t(X)=I+t(X-I)+4b_2t(1-t)(X-I)^2+O(\epsilon^3)\text{.}
\end{equation}
\end{prop}
\begin{pf}
Since $M_t(A,B)$ is a matrix mean if the generating $M(A,B)$ is a symmetric matrix mean therefore it has the following representation
\begin{equation}
M_t(A,B)=A^{1/2}f_t\left(A^{-1/2}BA^{-1/2}\right)A^{1/2}=Af_t\left(A^{-1}B\right)\text{,}
\end{equation}
where $f_t(X)$ is a normalized operator monotone function in $X$, therefore analytic on $(0,\infty)$, hence we have the second equality as well. Since it is generally analytic only on $(0,\infty)$, we expect \eqref{fexpansion} to be convergent only for $\left\|X-I\right\|<1$. By the above representation for $M_t(A,B)$ and the fact that $M_{1/2}(A,B)=M(A,B)$ by definition, it is enough to show that the expansion in the assertion holds for $f_t(X)$. In other words we have to consider the mean $M_t(A_0,B_0)$ of $A_0=I$ and an arbitrary $B_0=X$. We also have a natural expansion in the neighborhood of $I$ for the inverse function as
\begin{equation}\label{invexpansion}
X^{-1}=\sum_{k=0}^{\infty}(-1)^k(X-I)^k\text{,}
\end{equation}
which is convergent for $\left\|X-I\right\|<1$. Now in every step of the Weighted mean process we have to compute a symmetric mean of two matrices and by the assumption of the assertion we have $\left\|X-I\right\|\leq\epsilon<1$. Without loss of generality we may write $A_j$ and $B_j$ in the following forms
\begin{equation}
\begin{split}
A_j&=I+y_1^j(X-I)+y_2^j(X-I)^2+O(\epsilon^3)\\
B_j&=I+z_1^j(X-I)+z_2^j(X-I)^2+O(\epsilon^3)\text{.}
\end{split}
\end{equation}
Now we will make use of the above expansions to get an expansion for $M(A_j,B_j)$ up to the $O(\epsilon^3)$ term as follows
\begin{equation}
\begin{split}
&M(A_j,B_j)=A_jf(A_j^{-1}B_j)=\left(I+y_1^j(X-I)+y_2^j(X-I)^2+O(\epsilon^3)\right)\\
&f\left[A_j^{-1}\left(I+z_1^j(X-I)+z_2^j(X-I)^2+O(\epsilon^3)\right)\right]=\\
&=\left(I+y_1^j(X-I)+y_2^j(X-I)^2+O(\epsilon^3)\right)f\left[\left(I-y_1^j(X-I)+\right.\right.\\
&\left.\left.+((y_1^j)^2-y_2^j)(X-I)^2+O(\epsilon^3)\right)\left(I+z_1^j(X-I)+z_2^j(X-I)^2+O(\epsilon^3)\right)\right]
\end{split}
\end{equation}
where we have used \eqref{invexpansion} to express $A_j^{-1}$ and \eqref{fexpansion} to express $f(X)$ up to $O(\epsilon^3)$ terms. After some calculation and taking into account that the terms $(X-I)^k$ with $k\geq 3$ are of $O(\epsilon^3)$, we get that
\begin{equation}
M(A_j,B_j)=I+\frac{y_1^j+z_1^j}{2}(X-I)+\left[\frac{y_2^j+z_2^j}{2}+b_2(y_1^j-z_1^j)^2\right](X-I)^2+O(\epsilon^3)\text{.}
\end{equation}
Note that since $A_0=I$ and $B_0=X$ we have $y_2^0=0$ and $z_2^0=0$. Hence it can be easily proved by induction that the terms
\begin{equation}
\begin{split}
y_2^j=b_2p_j(y_1^0,z_1^0,t)\\
z_2^j=b_2q_j(y_1^0,z_1^0,t)\text{,}
\end{split}
\end{equation}
where $p_j$ and $q_j$ are functions which do not depend on $b_2$. Also since $f_t(X)$ is an analytic function due to Kubo-Ando theory, therefore the limits $p=\lim_{j\to\infty}p_j=\lim_{j\to\infty}q_j$ and this limit function is also independent of $b_2$.

Now since the weighted geometric mean
\begin{equation}
G_t(A,B)=A^{1/2}\left(A^{-1/2}BA^{-1/2}\right)^tA^{1/2}=A\left(A^{-1}B\right)^t
\end{equation}
is an affine mean, therefore the Weighted mean process gives back $G_t(A,B)$ for every $t\in[0,1]$ and $A,B\in P(r)$. In other words in this case if we expand the function $X^t$ into a Taylor series around $I$ we get
\begin{equation}
G_t(I,X)=X^t=I+t(X-I)-\frac{t(1-t)}{2}(X-I)^2+O(\epsilon^3)
\end{equation}
and this equation for $t=1/2$ gives that $b_2=1/8$. Since $G_t(A,B)$ is an affine mean and $p$ does not explicitly depend on $b_2$ we get $p=4t(1-t)$. Similar consideration can be applied in the case of the linear term $t(X-I)$.
$\square$
\end{pf}

The above proposition tells us that no matter how we choose the symmetric matrix mean $M(A,B)$, the series expansion of $M_t(A,B)$ around $I$ will have similar structure up to the $(X-I)^2$ term:
\begin{equation}
M_t(A,B)=(1-t)A+tB-4b_2t(1-t)A^{1/2}\left(A^{-1/2}BA^{-1/2}-I\right)^2A^{1/2}+\ldots
\end{equation}
Actually $0\leq b_2$ for every matrix mean since the corresponding operator monotone function is concave. It is not hard to prove using Proposition~\ref{arithmeticmatrixmean} that $0\leq b_2\leq 1/2$ for every matrix mean.

A remarkable consequence of the above series expansion is the following
\begin{thm}
The BMP procedure converges at least cubically if the matrices are sufficiently close to each other for all symmetric matrix means.
\end{thm}
\begin{pf}
This is so because we can mimic the proof of Theorem 3.2 in \cite{bini} which ensures the cubic convergence of the BMP process for the geometric mean. The only differences we need to take care of are in the series expansions. Namely that in property \textbf{C1} of Theorem 3.2 in \cite{bini} we will get
\begin{equation}
T_k:=\frac{1}{k}\sum_{j=1}^{k}E_j-\frac{2b_2}{k^2}\sum_{i,j=1}^{k}(E_i-E_j)^2\text{.}
\end{equation}
The other parts of the proof are just the same.
$\square$
\end{pf}

\bibliographystyle{elsarticle-num}

\end{document}